\documentclass[11pt]{article}
\usepackage{amsmath}
\usepackage{amsbsy}
\usepackage{amssymb}
\usepackage{latexsym}
\usepackage[lite]{amsrefs}
\def\r{{\mathcal R}}

\def\A{{\mathbb A}}
\def\C{{\mathbb C}}
\def\Q{{\mathbb Q}}
\def\N{{\mathbb N}}
\def\Z{{\mathbb Z}}

\def\aA{{\mathcal A}}
\def\iI{{\mathcal I}}
\def\jJ{{\mathcal J}}

\def\deg{{\rm deg}}
\numberwithin{equation}{section}
\usepackage[OT2,OT1]{fontenc}
\def\Sha{{\fontencoding{OT2}\selectfont SH}}

\begin{document}
\title{A Theorem on GL$(n)$ \`a l\`a Tchebotarev}
\author{
Dinakar Ramakrishnan\footnote{Supported by the Simons Foundation through a {\bf Simons Fellowship}, and partly by the NSF through an earlier grant DMS-1001916.}}
\date{}
\maketitle

\begin{center}
{\bf Abstract}
\end{center}

\begin{small}

\medskip

Let $K/F$ be a finite Galois extension of number fields. It is well known that the
Tchebotarev density theorem implies that an irreducible, finitely ramified $p$-adic representation
$\rho$ of the absolute Galois group of $K$ is determined (up to equivalence) by the
characteristic polynomials of Frobenius elements Fr$_v$ at any set of primes $v$ of $K$ of degree
$1$ over $F$. Here we prove an analogue for GL$(n)$, namely that a cuspidal
automorphic representation $\pi$ of GL$(n, \A_K)$ is determined up by the knowledge of
its local components at the primes of degree one over $F$. We prove in fact a stronger theorem, stimulated by a question of Michael Rapoport and Wei Zhang, relaxing to an extent the Galois hypothesis. The method uses, besides the Rankin-Selberg theory of L-functions
and the Luo-Rudnick-Sarnak bound for the Hecke roots of $\pi$, certain consequences
of class field theory via Galois cohomology. In an earlier paper (\cite{Ra2}) we obtained such a result up to twist equivalence for $K/F$ cyclic of prime degree by using basic Kummer theory. We make use of suitable solvable base changes $\pi_M$, relative to certain auxiliary succession of abelian extensions $E/F$, with $M$ being an abelian extension of the compositum $EK$, and deduce that $\pi_M \simeq \pi'_M$, and then descend this isomorphism to one over $K$. A key ingredient for progress here is the use of global Tate duality and a local-global result arising from class field theory. In fact we prove the main result for {\it isobaric} automorphic representations, which are analogues of {\it semisimple} Galois representations. In the last section we introduce a notion of {\it semi-temperedness}, which is much weaker than temperedness, but allows for the deduction of the main result without any hypothesis whatsoever on $K/F$.

\end{small}

\section{Introduction}

Let $F$ be a number field with adele ring $\A_F$. The object of this
article is to prove the following:

\medskip

\noindent{\bf Theorem A} \, \it Let $K/F$ be a finite extension of number fields with a filtration
$$
F=K_0 \subset K_1 \subset \dots \subset K_m =K,
$$
where each $K_j$ is a field normal over $K_{j-1}$ ($\forall j\in\{1,\dots, m\}$). Let
$\Sigma^1(K/F)$ be the set
of primes of $K$ which are of degree $1$ over $F$. Fix any $n\geq 1$, and consider
isobaric automorphic representations $\pi =\otimes_v \pi_v$, $\pi'=\otimes_v \pi'_v$ of GL$_n(\A_K)$. Suppose
that $\pi_v \simeq \pi'_v$ for all but a finite number of $v$ in $\Sigma^1(K/F)$.
Then $\pi \simeq \pi'$.
\rm

\medskip

Note that $K/F$ is not assumed to be Galois or solvable. This result is analogous to, and inspired by, the well known consequence of the Tchebotarev density theorem that a semisimple, finitely ramified, $n$-dimensional $\overline \Q_p$-representation $\rho$ of the absolute Galois group of $K$ is determined, up to isomorphism, by the collection of its restrictions $\rho_v$ to the (decomposition groups of) primes $v$ of degree $1$ over $F$. Note that our main result applies in particular to cuspidal automorphic representations $\pi$ of GL$_n(\A_K)$, which are the building blocks of isobaric representations (cf. \cite{La}, \cite{JS}), and are expected, when algebraic, to be associated to irreducible, potentially semistable $p$-adic Galois representations of dimension $n$ which are unramified outside a finite set of primes. Even if one is interested only in the case when $\pi, \pi'$ are cuspidal, the fact that they may become Eisensteinian when we apply base change to suitable solvable extensions necessitates working in a larger framework. On the other hand, it is essential to stick to the isobaric setup of Langlands for otherwise one can produce (non-cuspidal) counterexamples using non-isobaric forms.

It may be useful to note that our result is non-obvious even in the case of twists, i.e., when $\pi'$ is of the form $\pi\otimes\chi$, for a character $\chi$ of the idele class group of $K$.

As it will be apparent, Theorem A is a strengthening of the celebrated, and oft-used {\it strong multiplicity one theorem} due to Jacquet, Shalika, and Piatetski-Shapiro, which deals with $\pi, \pi'$ agreeing outside a {\it finite} number of places. One can also effective bounds, and this was done in \cite{Mo}, which in fact began the analytic approach to this problem; a better bound has recently been established in \cite{B}.

When $n=1$, automorphic forms are just Grossencharacters $\chi$, and a theorem of Hecke asserts that $\chi$ is even determined by the knowledge of $\chi_v$ at a set of primes $v$ of density $> 1/2$. An analogue of this is known for GL$(2)$ (\cite{Ra1}), where a cusp form is determined by its components at any set of primes of density $> 7/8$; see also \cite{Wa}. We are far from any such result for GL$(n)$ for $n \geq 3$, as functoriality is not known for the adjoint $L$-function. There is an elegant Galois analogue for $\ell$-adic representations of any dimension due to Rajan (\cite{Rj1}).

A milder version of Theorem A was established for GL$(n)$ by the author in 2010 (\cite{Ra2}), which dealt only with cyclic extensions $K/F$ and also determined $\pi$ only up to twist equivalence (unless $[K:F]=2$). We refer to that article for the basic analytic argument, which is also of utility here. The difficulty caused by the lack of knowledge of the general Riemann Hypothesis for unitary cusp forms $\pi$ on GL$(n)/K$ is to a degree alleviated by the bound of Luo, Rudnick and Sarnak (\cite{LRS}). Its use brings into focus the the set $\Sigma=\Sigma(\pi,\pi')$ of finite places $v$ of degree $j$ over $F$, with $2 \leq j \leq n^2+1$, where $\pi_v \not\simeq \pi'_v$. The main point of this paper is its resolution, and the proof is algebraic, using the global Tate duality of Class field theory, and a local-global result.

Our basic idea is to move to a bigger Galois extension $L$ of $F$ containing $K$ with $L/K$ solvable, and with $[L:F]$ being divisible by the same primes as $[K:F]$, such that the divisors $\tilde v$ in $L$ of those in $\Sigma$ have higher degree over $F$. Roughly speaking, we use an inductive argument and at each stage the relevant $\Sigma$ is partitioned into a finite union of subsets $\Sigma_j$, with all but finitely many places in each $\Sigma_j$ acquiring a higher degree in a suitable Galois extension $L_j/F$ containing $K$ (and solvable over $K$). There will of course be new places $u$ of $L_j$ with low degree, but they will be arising from degree one places in $K$ and so we would know that the base changes $\pi_{L_j}$ and $\pi'_{L_j}$ agree at such $u$.

It may be instructive to look at the simplest case, namely when $K/F$ is a cyclic extension of prime degree $p$, with $F$ containing a primitive $p^2$-th root of unity, and such that $p < \frac{n^2+1}{2} \leq p^2$. Here, after writing $K$ as $F[\alpha^{1/p})]$ for some $\alpha\in F^\ast-{F^\ast}^p$, we may consider $L=F[\alpha^{1/p^2}]$, in which the places of $F$ which become inert in $K$ stay inert and hence acquire degree $p^2$ (over $F$), allowing us to conclude that the base changes (cf. \cite{AC}) $\pi_L$, $\pi'_L$ to GL$(n)/L$ are isomorphic. Suppose for simplicity of exposition here that $\pi, \pi'$ are cuspidal. Then $\pi'\simeq \pi\otimes\delta$, for a character $\delta$ of $K$ becoming trivial on $L$ (when pulled back by norm). If $\delta$ is trivial, there is nothing to prove, so we may take $\delta$ to cut out $L/K$. Besides, we may replace $L$ by $L'=F[(\alpha\beta)^{1/p^2}]$ with $\beta\in {F^\ast}^p-{F^\ast}^{p^2}$, which furnishes a similar isomorphism $\pi'\simeq\pi\otimes\delta'$. Putting them together, we may assume that $\pi$ admits a non-trivial self-twist under $\lambda:=\delta'\delta^{-1}$, implying in particular that $p \mid n$, in turn forcing $n=p$ since $2p^2\geq n^2+1$. By varying $\beta$, we get many different such self-twists, one for each element of the image of $F^\ast$ in $K^\ast/{K^\ast}^p$, which is huge. To get a contradiction, we check (see Lemma B in section 4) that the number of such self-twists for any cusp form on GL$(n)/K$ is bounded above by $n^2$.

The general case is more involved, and the results and methods of class field theory play a key role. We repeatedly make use of auxiliary abelian extensions and auxiliary places, and we also appeal to the solvable base change for GL$(n)$ (\cite{AC}), to deduce that base changes of $\pi$ and $\pi'$ to suitable solvable Galois extensions $R$ of $K$, with $R/F$ Galois, are isomorphic. Then we vary $R$ suitably by making an auxiliary finite set of places split completely in $R$, allowing us to eventually descend the isomorphism to on over $K$.

If $F$ were to contain sufficiently many roots of unity, then by using Kummer theory, the representations $\pi, \pi'$ satisfying the hypothesis of Theorem A can be shown to be isomorphic by a finer version of the argument given above for $p$-extensions. When the roots of unity are not present, the most one can show by base changing to appropriate cyclotomic fields is that $\pi, \pi'$ are twist equivalent. (This is what we achieved in \cite{Ra2} for $K/F$ cyclic, but the proof here is different, and self-contained.) To do better and eliminate this twisting ambiguity, we have to work without worrying about the presence of roots of unity. Then the difficulty of writing certain Galois characters as $p$-th powers, and this leads us to contend with some obstructions in Galois cohomology (in degree $2$). For every $p$ dividing $\vert {\rm Gal}(K/F)\vert$, we appeal to the description of $H^2(F,\Z/p^m)$ in \cite{ArT} via the Tate duality concerning the global-to-local kernel. ($H^2(F,\Z/p^m)$ is the $p^m$-part of the Brauer group of $F$ when the $p^m$-th roots of unity are in $F$, but not otherwise.) The obstructions appear over intermediate fields $K'$ with $[K:K']=p$, which we trivialize over auxiliary extensions $M'$ in $M=KE$, with $[M:M']=p$ and $E/F$ abelian. This way we avoid assuming that $F$ contains any root of unity, and the choice of $E/F$ is sufficiently flexible to force, in addition, any given finite set $T$ of good places $v$ in $K$ to split completely in $M$. This allows us to descend later to $K$ the isomorphism of (the base changes of) $\pi$ and $\pi'$ over $M$. It is important to note that at no point do we assume that $K/F$ is solvable.

We also have the following consequence of Theorem A (when combined with solvable base change cf. \cite{AC} for GL$(n)$):

\medskip

\noindent{\bf Corollary B} \, \it Let $K/F$ be a finite extension of number fields, whose normal closure $\tilde K$ over $F$ is solvable over $K$, but not necessarily over $F$. Then, if $\pi, \pi'$ are cuspidal automorphic representations of GL$_n(\A_K)$ have isomorphic local components at all but a finite number of primes of degree one over $F$, then their base changes $\pi_{\tilde K}$ and $\pi'_{\tilde K}$ (to GL$(n)/\tilde K$) are isomorphic.\rm

\medskip

There are two kinds of potential uses for Theorem A. One is in establishing functoriality from a reductive group $G/F$ (such as a unitary group) which becomes isomorphic to GL$(n)$ an extension field $K$, providing sufficiency of checking at the primes of $K$ which are of degree one over $F$. (When $G$ is defined by a division algebra of dimension $n^2$ over a CM field $K$ with an involution of the second kind, one needs to first make a quadratic base change to $K$, where $G_K=D^\times$, and then make a further base change of degree $n$ to get to GL$(n)$.) Another potential utility is in the simplification one has in the local moduli problem when the place of interest is of degree one over $\Q$. We refer to the articles \cite{FMW}, \cite{GWa}, \cite{Zh}, and the forthcoming paper of Michael Rapoport and Wei Zhang.

\medskip

Finally, let us call an isobaric automorphic representation $\pi$ of GL$_n(\A_K)$ {\it semi-tempered} iff there exists $\delta >0$ such that for all but a finite number finite places $v$ (with norm $N(v)$) where $\pi$ is unramified, the associated Langlands class $\{\alpha_{1,v}, \dots, \alpha_{n,v}\}$ satisfies $\vert \alpha_{j,v}\vert \, \leq \, CN(v)^{1/4-\delta}$.

\medskip

\noindent{\bf Theorem C} \, \it Let $K$ be a number field and $\pi, \pi'$ isobaric automorphic representations of GL$_n(\A_K)$.
\begin{enumerate}
\item[(a)]Suppose $\pi, \pi'$ are semi-tempered and satisfy $\pi_v \simeq \pi'_v$ at all but finitely many finite places $v$ of degree one (over $\Q$). Then $\pi$ and $\pi'$ are globally isomorphic.
\item[(b)]Suppose there is an isobaric automorphic representation $\Pi$ of GL$_{n^2-1}(\A_K)$ such that the standard $L$-function $L^S(s, \Pi)$ equals $L^S(s,\pi, {\rm Ad})$ (for some finite set $S$ of places), where ${\rm Ad}$ denotes the adjoint representation of GL$_n(\C)$. Then $\pi$ is semi-tempered.
\end{enumerate}
\rm

\medskip

A proof of this result is in the last section. Part (a) follows from the slight strengthening of section 2 of \cite{Ra2} facilitted by the hypothesis of being semi-tempered. The basic analytic setup used there was considered already by Hecke for $n=1$, and for general $n$ by C.~Moreno (\cite{Mo}), made more effective in \cite{B}. Part (b) is a consequence of the bound of Luo, Rudnick and Sarnak.It should be remarked tat the existence of the Adjoint transfer $\pi \mapsto {\r Ad}(\pi)$ is one of the most pressing problems in automorphic forms, which is known for $n=2$ by Gelbart and Jacquet and for $n=3$ and $\pi$ essentially selfdual by the symmetric power lifting results of for GL$(2)$ by Kim and Shahidi. It is open already for $n=3$ and $\pi$ non-selfdual.

\medskip

We thank the participants, in particular H.~ Jacquet and F.~Shahidi,  at the memorial conference some years back for Joseph Shalika at the Johns Hopkins University, where we explained the cyclic case of our result, and we thank P.~Michel for invitation to give a lecture at EPFL in Lausanne, where we discussed a preliminary version of this paper. Essentially the final version of the paper was presented at the International Colloquium on Automorphic Forms and L-functions in January 2012 at the Tata Institute of Fundamental Research in Mumbai, where we had interesting conversations with many participants including D.~Goldfeld, C.S.~Rajan and D.~Rohrlich. The writing of the manuscript has benefited from the comments of all these mathematicians. We acknowledge partial support from the NSF through the grant DMS-1001916. The final version, with an improved result, was done while the author was a Simons Fellow, visiting Princeton University, whose hospitality he thanks, especially Chris Skinner, as well as encouraging communication since then with Wei Zhang.

\bigskip

\section{The first reduction and twists}

\medskip

For any global field $F$ with ad\`ele ring $\A_F$, let $\Sigma_F$
denote the set of all
places of $F$.  If $v\in \Sigma_F$ is finite, let $q_v$ denote the cardinality
of the residue field at $v$.  For $n \geq 1$, let $\aA_0(n,F)$, resp. $\aA_u(n,F)$, denote the set of
isomorphism classes irreducible, cuspidal, resp. isobaric sum of unitary cuspidal, automorphic
representations $\pi = \otimes_v \,
\pi_v$ of GL$(n,\A_F)$ with $\pi_v$ is unramified at all $v$ prime to the conductor of $\pi$.
We refer to the first three sections of \cite{Ra2} for the basic facts we need to establish the Theorem.

For every finite extension $L$ of $F$, and for all $j \geq 1$, put
$$
\Sigma_{L/F}^j: \, = \, \{v\in \Sigma_F \, \vert \, \deg_F(v)=j\},
$$
where $\deg_F$ denotes the degree over $F$.

\medskip

\noindent{\bf Proposition 1.1} \, \it Let $K/F$, $\pi, \pi'$ be as in Theorem A. Then for every finite subset $T$ of $\Sigma^1_{K/F}$ prime to the conductor of $\pi, \pi'$ and the absolute discriminant of $K$, there exists a finite extension $M/K$ with a filtration
$$
K=M_0 \subset M_1 \subset \dots \subset M_r=M
$$
with each $M_j/M_{j-1}$ cyclic, such that
\begin{enumerate}
\item[(a)]The base changes $\pi_M, \pi'_M$ of $\pi, \pi'$ respectively (to $M$) are isomorphic;
\item[(b)]The places in $T$ split completely in $M$.
\end{enumerate}
\rm

\medskip

For the existence of $\pi_M, \pi'_M$ for such an extension $M/K$ (which need not be Galois, but solvable), we refer to \cite{AC}.

\medskip

\noindent{\bf Claim 1.2} \, \it Proposition 1.1 \, $\implies$ \, Theorem A.
\rm

\medskip

This would be easy to prove if in Prop. 1.1 the degrees of the places of $T$ could be unrestricted.

\medskip

{\it Proof of Claim} \, Since the central characters $\omega, \omega'$ of $\pi, \pi'$ respectively agree at a set of primes of density one, we know by Hecke that necessarily, $\omega=\omega'$.

We may write $\pi$ uniquely as an isobaric sum
$$
\pi \, \simeq \, \pi_1 \boxplus \dots \boxplus \pi_r
$$
which each $\pi_i$ a cusp form on GL$(n_i)/F$, with $\sum_i n_i = n$, allowing $\pi_i \simeq \pi_j$ for $i=j$. Similarly,
$$
\pi' \, \simeq \, \pi'_1 \boxplus \dots \boxplus \pi'_s
$$
which each $\pi'_j$ a cusp form on GL$(n'_j)/F$, with $\sum_j n'_j = n$.

We will use three stages of nested induction. To begin, since we know the claim for $n=1$, we will let $n>1$ and assume by induction that the assertion holds for all $k < n$.

Put
$$
m \, = \, [M:K].
$$
If $m=1$, there is nothing to prove. So let $m>1$ and assume that the Claim holds for all $m' < m$.
We may refine the filtration $\{M_j\}$ such that each $M_j/M_{j-1}$ is cyclic of prime degree, with $M_0=K$. Now consider the primes $u$ of degree $1$ of $M_1$ over $F$. Then the primes $v$ of $K$ below them have degree one over $F$, and by hypothesis, $\pi_v \simeq \pi'_v$ for every such that, and for each $u \vert v$, the basic property of base change implies that $\pi_{M_1,u} \simeq \pi'_{M_1,u}$. So we may apply Prop. 1 with $M_1$ in the place of $K$, whose conclusion implies by induction, since $[M:M_1]$ is smaller than $m=[M:K]$, that
$$
\pi_{M_1} \, \simeq \, \pi'_{M_1},
$$
with $[M_1:K]=\ell$, a prime.
Consequently,
if $M_1$ is cut out by a character $\delta$ of $K$, we have
$$
\boxplus_{k=0}^\ell \boxplus_{i=1}^r \, \pi_i\otimes\delta^k \, \, \simeq \, \,
\boxplus_{k=0}^\ell \boxplus_{j=1}^s \, \pi'_j\otimes\delta^k.
$$
This is because
$$
I_{M_1}^K(\underline{1}) \, \simeq \, \boxplus_{k=0}^{\ell-1} \, \delta^k,
$$
where the left hand side is the automorphic induction of the trivial representation of GL$(1)/M_1$ to GL$(\ell)/K$, and
$$
L(s, \eta \times I_{M_1}^K(\underline{1})) \, = \, L(s, \boxplus_{k=0}^{\ell-1} \, \eta\otimes\delta^k),
$$
for any isobaric automorphic representation $\eta$ of GL$(n)/K$.

It follows, by the uniqueness of the isobaric sum decomposition (\cite{JS})
that
$$r \, = \, s,
$$
and for all $(i,k)$ with $i \leq r$, $k \leq \ell-1$, there is an index $j=j(i)\leq r$ and an integer $a \in \{0, \dots, \ell-1\}$ such that
$$
\pi_i\otimes\delta^k \, \simeq \, \pi'_j\otimes\delta^a.
$$
Let us rewrite this as
$$
\pi'_j \, \simeq \, \pi_i \otimes\delta^b, \, \, \, \forall \, i\leq r,\leqno(\ast)
$$
with $j=j(i)$ and $b=b(i) \in \{0, \dots, \ell-1\}$.

Suppose $r=1$. Then $\pi, \pi'$ are cuspidal with
$$
\pi' \, \simeq \, \pi\otimes \delta^b.
$$
If $b=0$, we are done. So let $b \in\{1, \dots, \ell-1\}$. Note that the extension $M_1/K$ depends on $T$, and we may choose a different finite subset $\tilde T$, say of $\Sigma^1_{K/F}$ which consists of primes which are inert in $M_1$. Then the extension $\tilde M_1$ attached to $\tilde T$ will be a disjoint cyclic extension of $K$, cut out by a character $\tilde \delta$. Then we will get
$$
\pi \otimes \delta^b \, \simeq \, \pi' \, \simeq \, \pi\otimes\tilde\delta^{\tilde b},
$$
implying that $\pi$ admits a self-twist by $\mu=\delta^{-b}\tilde\delta^{\tilde b}$.
There are infinitely many choices for $T$ and we may choose them so that the corresponding cyclic extensions of $K$ are all disjoint. This will result, if $\pi'$ were not isomorphic to $\pi$, an infinite number of self-twists of $\pi$, which is impossible. (There can be at most $n^2$ possibilities since for every such twisting character $\nu$, $L(s, \pi\times \overline\pi\otimes\nu^{-1})$ must have a pole at $s=1$.
It follows that $\pi$ and $\pi'$ must be isomorphic over $K$ when $r=1$.

So we may take $r>1$ and assume by induction the assertion for the isobaric representations of length $< r$. So it suffices to show that $\pi_i \simeq \pi'_j$ for some $i$, and corresponding $j=j(i)$ (see $(\ast)$. The same argum,ent as in the cuspidal case shows, by choosing infinitely many disjoint cyclic extensions $M_1/K$ corresponding to different choices of $T$, that (by the pigeon hole principle) for at least one $i$, call it $i_0$, there is a $j_0$ such that $\pi_i$ is isomorphic to $\pi'_j$ for some $j$. Then we only have to prove that
$$
\boxplus_{i\ne i_0} \, \pi_i \, \simeq \, \boxplus_{j\ne j_0} \, \pi'_j,
$$
which follows by induction as the length is now $r-1$.

Done.

\bigskip

\section{The Second Reduction}

\medskip

For any number field $L$ containing $F$, and for each $j \geq 1$, let $\Sigma_{L/F}^j$ denote the set of primes of $L$ which are of degree $j$ over $F$.

\medskip

\noindent{\bf Proposition 2.1} \, \it Let $K/F$, $\pi, \pi'$ be as in Theorem A. Then for every finite set $T$ of finite places of $K$ prime to the conductor of $\pi, \pi'$ and the discriminant of $K$, there exists a finite extension $M/K$ with a filtration $K=M_0 \subset M_1 \subset \dots \subset M_r=M$ with each $M_j/M_{j-1}$ cyclic, such that
\begin{enumerate}
\item[(a')]For every $j \leq[(n^2+1)/2]$, for all primes $v$ in $\Sigma_{K/F}^j$ outside a finite set, and for every prime $u$ of $M$ above $v$, either $u$ has degree $> [(n^2+1)/2]$ over $F$ or $\pi_{M,u}\simeq \pi'_{M,u}$;
\item[(b')]The places in $T$ split completely in $M$.
\end{enumerate}
\rm

\medskip

\noindent{\bf Claim 2.2} \, \it Proposition 2.1 \, $\implies$ \, Theorem A.
\rm

\medskip

{\it Proof of Claim} \, Thanks to Claim 1.2, it suffices to show that Prop. 2.1 implies Prop. 1.1. Since part (b') of Prop. 2.1 is identical to part (b) of Prop. 1.1, we need only the check that part (a') implies part (a) of Prop. 1.1. It is in effect a consequence of the bound on the Hecke eigenvalues due to \cite{LRS}, once we admit the existence of $M$ as in Prop. 2.1. Indeed, it suffices to show, thanks to Prop. 2.1 of \cite{Ra2} (with $M$ playing the role of $F$ there), that for {\it all but finitely many primes} $u$ of $M$ {\it of degree $\leq [(n^2+1)/2]$}, we have
$$
\pi_{M,u} \, \simeq \, \pi'_{M.u}.\leqno(\ast)
$$
If $v$ is a place of $K$ below such a $u$, its degree (over $F$) is necessarily at most $[(n^2+1)/2]$. If $v$ has degree $1$ over $F$, then by the hypothesis of Theorem A, $\pi_v$ and $\pi'_v$ are isomorphic (outside a finite number of exceptions), yielding $(\ast)$ for any $u$ of $M$ above a degree $1$ prime of $K$. ($u$ need not have degree one over $F$.) So consider when $v$ has degree between $2$ and $[(n^2+1)/2]$. Then by Prop. 2.1, the degree of the place $u$ (of $M$) above $v$ is $> [(n^2+1)/2]$, so irrelevant for $(\ast)$. Hence $(\ast)$ holds for all but a finite number of $u$ of degree $\leq [(n^2+1)/2]$. By Prop. 2.1 of \cite{Ra2}, $\pi_M$ and $\pi'_M$ are isomorphic, proving Prop. 1.1 (of this paper).

Done.

\bigskip

\section{The width and the third reduction}

\medskip

The base number field $F$ will be fixed throughout.
Put
$$
h\, = h(\pi,\pi'; K/F) \, = \, \max\{j  \leq [\frac{n^2+1}{2}]+1 \, \, \, \, : \, \, \, \, \forall i < j, \pi_v \simeq \pi'_v, \, \forall \, v\in \Sigma^i_{K/F}-Z_i, \, \vert Z_i\vert<\infty\}.\leqno(7.2)
$$

The {\it width} of $(\pi, \pi'; K/F)$ is defined as follows:
$$
w = w(\pi, \pi'; K/F)\, = \, [\frac{n^2+1}{2}]+1-h \, \, \in \, \Z.\leqno(7.3)
$$
Since $1\leq h \leq [\frac{n^2+1}{2}]+1$,
$$
0 \, \leq \, w \, \leq \, [\frac{n^2+1}{2}].
$$
Note that when we make a base change to a finite extension $\tilde K/K$, which is filtered by cyclic extensions (so that the base change exists by \cite{AC}), the width does not increase, i.e.,
$$
w(\pi_{\tilde K}, \pi'_{\tilde K}; \tilde K/F) \, \leq \, w(\pi, \pi'; K/F),
$$
but it could potentially remain the same.

\medskip

\noindent{\bf Proposition 3.1} \, \it Let $K/F$, $\pi, \pi'$ be as in Theorem A. Suppose $w>0$. Then for every finite set $T$ of finite places of $K$ prime to the conductor of $\pi, \pi'$ and the discriminant of $K$, there exists a finite extension $N/K$ filtered as $K=N_0 \subset N_1 \subset \dots \subset N_r=N$ with each $N_j/N_{j-1}$ cyclic, such that
\begin{enumerate}
\item[(a)]$w(\pi_N, \pi'_N; N/F) \, < \, w(\pi, \pi'; K/F)$;
\item[(b)]The places in $T$ split completely in $N$.
\end{enumerate}
\rm

\medskip

\noindent{\bf Claim 3.2} \, \it Proposition 3.1 \, $\implies$ \, Theorem A.
\rm

\medskip

{\it Proof of Claim} \, In view of Claim 2.2, it suffices to deduce Prop. 2.1. For this we will use induction on the width. If $w=0$, then we may take $N=K$ and the assertion follows.
So we may let $w >0$ and assume that the assertion holds for smaller widths. Let $N$ be as in Prop. 3.1. Since $w(\pi_N, \pi'_N; N/F)$ is (strictly) smaller than $w$, Proposition 2.1 holds by induction for $N$ in the place of $K$ there, such that if $\tilde T$ is the set of places of $N$ above $T$, then it splits completely in $M$. Now since $M/N$ and $N/K$ are both filtered by successive cyclic extensions, part (a') of Prop. 2.1 holds for $(\pi, \pi')$ over $K$ as well. Moreover, since by construction, every place in $T$ splits completely in $N$, and since every place above it (in $\tilde T$) splits completely in $M$, part (b') of Prop. 2.1 also holds.

Done.

\section{Trivialization of certain torsion classes in Galois cohomology}

\medskip

The following Lemma will play a key role for us, and it is needed partly to avoid assuming that $F$ contains sufficiently many roots of unity, and partly to have flexibility in our choice relative to the auxiliary finite set $T$ of places in $K$ of degree one over $F$.

\medskip

\noindent{\bf Lemma 4.1} \, \it
Let $K/F$ be a finite Galois extension of number fields, with $p$ a prime divisor of the order of Gal$(K/F)$. Let $\iI$ be the set of intermediate fields $K'$ in the extension $K/F$ such that $[K:K']=p$. For each $K'\in \iI$, fix a class $\beta(K')$ in $H^2(K',\Z/p)$ whose restriction to $K$ is trivial (in $H^2(K, \Z/p)$). Fix also an auxiliary finite set $T_0$ of finite places of $F$ which are prime to $p$ and split completely in $K$, and let $T$ be the places of $K$ above $T_0$.Then there exists a cyclic extension $E/F$ such that
\begin{enumerate}
\item[(a)]$E/F$ is linearly disjoint from $K/F$, and is of degree $p^r$ or $2p^r$,
with $r$ being independent of $T$;
\item[(b)]Every place in $T$ splits completely in $KE$;
\item[(c)]For every $K'\in \iI$, and for every subfield $M'$ of $KE$ containing
$K'$ with $[KE:M']=p$, the restriction of $\beta(K')$ in $H^2(M',\Z/p)$ is trivial.
\end{enumerate}
\rm

\medskip

{\it Proof}. \, For $0\leq i\leq 3$, consider the map
$$
\alpha_i=(\alpha_{i,u}): H^i(K', \Z/p) \, \rightarrow \, \prod_u \,H^i(K'_u,\Z/p),
$$
where $u$ runs over all the places of $K'$, $K'_u$ denotes the local completion of $K'$ at $u$, and $\alpha_u$ the restriction at $u$. It is known (see \cite{Mi}, chapter 1) that for any $\beta\in H^2(K',\Z/p)$, there is a finite set $X(K')$ of places of $K'$ such that $\alpha_{i,u}(\beta)$ is zero at every place $u$ outside $X(K')$, which is seen by noting that $\beta$ must be in the image of $H^2({\rm Gal}(L/K'),\Z/p)$ for a finite Galois extension $L/K'$. Moreover, the kernel \Sha$^i(K',\Z/p)$ of $\alpha_i$ is, by Tate, in duality with \Sha$^{3-i}(K',\mu_{p})$
(\cite{Mi}), where $\mu_{p}$ denotes the Galois module of $p$-th roots of unity; $\Z/p$ is as usual the trivial Galois module. By Artin-Tate \cite{ArT}, \Sha$^1(K',\mu_{m})$, and hence Sha$^2(K',\Z/m)$, is either trivial or of order $2$. In fact, for $m=p$, this kernel is trivial.

\medskip

\noindent{\bf Sublemma 4.2} \, \it For every $u\in X(K')$, the restriction of $\beta_u:=\alpha_{2,u}(\beta(K'))$ becomes trivial over the unique unramified $p$-extension $L_1(u)$, say, of $K'_u$.\rm

\medskip

{\bf Proof}. \, We may (and we will) assume from here on that $\beta_u$ is non-trivial on $K'_u$ for any $u\in X(K')$. Indeed, if we put $k=K'_u[\mu_{p}]$, then the restriction $\beta_{u,k}$ of $\beta_u$ to $k$ lies in $H^2(k,\mu_{p})=Br_k[p]$, the $p$-torsion subgroup ($\simeq \Z/p$) of the Brauer group of $k$. It is well known (cf. \cite{Se}) that any class in $Br_k[p]$ becomes trivial over the unramified $p$-extension $k'$ of $k$. Let $\delta$ denote the unramified character of ${K'_u}^\ast$ of order $p$. The pull-back by norm of $\delta$ to $k$ is evidently non-trivial and cuts out the extension $k'/k$. Let $\beta_{u,1}$ be the restriction of $\beta_u$ to $L_1(u)$. Then its restriction to $k'$ is trivial, and since the composition of restriction followed by norm is multiplication by the degree, we see that $[k':L_1(u)]\beta_{u,1}=0$. Since$[k':L_1(u)]=[k:K'_u]$ divides $p-1$, and since $\beta_{u,1}$ is killed by $p$, we must have $\beta_{u,1}=0$, proving the assertion of the Sublemma.

\medskip

{\bf Proof of Lemma 4.1} ({\it contd.}) \, Now let $X_0$ denote the finite set of places of $F$ above which lie all the places of $X(K')$ for all $K'\in \iI$. Since $\beta_u$ is non-trivial on $K'_u$ for any $u\in X(K')$, $u$ cannot split in $K$ (as $\beta$ becomes trivial over $K$). Let $v$ be the unique place over $u$; $K_v/K'_u$ may or may not be ramified. In any case choose, for each $u_0\in X_0$, a finite, cyclic, unramified extension $E(u_0)$ of $F_{u_0}$, say of degree $p^r$, for large enough $r>0$, such that for every $u$ above $u_0$ lying in some $X(K')$, the compositum $K_vE(u_0)$  contains $L_2(u)$, the unramified (cyclic) $p^2$-extension of $K'_u$. (For each $u_0$, there is an $r>0$ which works, and since there are only a finite number of $u_0$'s in $X_0$, we may choose an $r$ which works for all of them.) In particular, the restriction of $\beta_u$ to the unramified extension ${K'_u}E(u_0)$ is trivial ($\forall K'\in \iI$). Now, by appealing to the Grunewald-Wang theorem (\cite{ArT}), we may choose a global cyclic extension $E/F$ of degree $tp^r$, with $t\in \{1,2\}$, such that
\begin{enumerate}
\item[(i)] the local extension of $E/F$ at any divisor of $u_0$ is $E(u_0)/F_{u_0}$,
\item[(ii)] $E/F$ and $K/F$ are linearly disjoint from each other, and
\item[(iii)] every place in $T_0$ splits completely in $E$.
\end{enumerate}
Then $KE$ is Galois over $F$, and hence over $K'$ for every $K'\in iI$. Also, $KE$ is abelian over $K'$ with Galois group $(\Z/p)\times(\Z/tp^r)^2$. By construction, $u\in X(K')$ is unramified in $K'E$, splitting into a product of places $u_1, \dots, u_m$ (of places in $K'E$) such that the restriction of $\beta_u$ to each $(K'E)_{u_j}$ is trivial, implying that the global class $\beta$ restricts to $0$ in $H^2(K'E, \Z/p)$, and hence in $H^2(K'E, \Z/p)$.
We have to prove furthermore that $\beta$ has trivial restriction to any intermediate field $M'$ in $KE/K'$ with $[KE:M']=p$, not just to $K'E$. Fix such an $M'$, which will be normal over $K'$ since $KE/K'$ is abelian. Consider any place $\tilde u$ of $M'$ above $u$ in $X(K')$. As noted above, since $\beta$ becomes trivial upon restriction to $K$, and since $\beta_u$ is non-trivial on $K'_u$, $u$ cannot split in $K$, and we write $v$ for the unique place of $K$ above $u$, and $\tilde v$ a place of $KE$ above $v$ such that $\tilde u$ lies below it in $M'$. By construction, $(KE)_{\tilde v}$ contains $L_2(u)$. It follows that, since $[KE:M']=p$, the local field $M'_{\tilde u}$ must at least contain $L_1(u)$, over which the local class $\beta_u$ becomes trivial (by Sublemma 4.2). So $\beta_u$ will become trivial over $M'_{\tilde u}$. This holds for every $u\in X(K')$, and so the global class $\beta$ restricts to zero over $M'$, as asserted.

Done.

\bigskip

\section{The main result in the Galois case}

\medskip

In this section we will prove the main result, which follows from Prop. 3.1, in this section for the critical case when $K/F$ is Galois. In fact we will show the following

\medskip

\noindent{\bf Proposition 5.1} \, \it Let $K/F$ be a finite Galois extension, and $T$ a finite subset of $\Sigma_{K/F}^1$.  Fix an integer $j$ in $\{2, \dots, [(n^2+1)/2]\}$. Then there exists a finite solvable extension $M/K$, realized as a successive extension of cyclic extensions of prime degree, such that
\begin{enumerate}
\item[(a)]The places in $T$ split completely in $M$;
\item[(b)]For all but a finite number of $v \in \Sigma_{K/F}^j$ with , we have, for every place $u$ of $M$ dividing $v$,
$$
\deg_F(u) \, > \, \deg_F(v).
$$
\end{enumerate}
\rm

\medskip

\noindent{\bf Claim 5.2} \, \it Prop. 5.1 \, $\implies$ \, Prop. 3.1 for $K/F$ Galois.\rm

\medskip

{\it Proof}. \, Let $K/F$, $\pi, \pi'$ be as in Theorem A. By hypothesis, we have $\pi_v \simeq \pi'_v$, for all but a finite number of places $v$ of $K$ of degree $1$ over $F$. By the basic property of base change, we then have
$$
\pi_{K,u} \, \simeq \, \pi'_{K,u}, \, \, \forall u\mid v, \, \, \forall \, v \in \Sigma_{K/F}^1.
$$
(Of course $u \in \Sigma_M$ need not have degree one over $F$.)

In view of part (a) Prop. 5.1, which we apply with $j=h$ \, ($=h(\pi, \pi'; K/F)$), and the definition of the width $w$ as $[(n^2+1)/2]+1-h$, we see that
$$
w(\pi_M, \pi'_M: M/F) \, < \, w(\pi, \pi'; K/F),
$$
yielding part (a) of Prop. 3.1. Part (b) is the same in both Propositions.

Done.

\medskip

{\it Proof of Prop. 5.1} \,

Let $v_0$ be a place in $\Sigma_{K/F}$, unramified over $F$, with degree equal to $j$, which is $\leq (n^2+1)/2$; so $N_{K/F}(v_0)=p^j$, where $p$ is a rational prime. (The only primes which intervene are those dividing $[K:F]$.) The decomposition group of $v_0$ over $p$ is a subgroup $H$ of Gal$(K/F)$, which, since $v_0$ is unramified over $F$, is cyclic of order $p^j$, and hence contains a unique subgroup $C$ of order $p$. Viewing $C$ as a subgroup of the global Galois group, we see that it defines a subfield $K'$ of $K$ with $[K:K']=p$.

Let $\iI$ denote (as before) the finite collection of intermediate fields $K'$ of $K/F$ such that $[K:K']=p$. For each such $\iI$, denote by $\varphi(K')$ the character of order $p$ of the absolute Galois group $\Gamma_{K'}$ of $K'$ cutting out the cyclic $p$-extension $K/K'$. The surjective $p$-power map $z \mapsto z^p$ on $\C^\ast$ gives a short exact sequence of {\it trivial} $\Gamma_{K'}$-modules
$$
0\to \Z/p \to \C^\ast\to \C^\ast \to 1.
$$
(We are looking at the trivial Galois action since we are interested in ${\rm Hom}(\Gamma_{K'},\C^\ast)$.) The associated long exact sequence in Galois cohomology yields
$$
H^1(\Gamma_{K'},\C^\ast) \to H^1(\Gamma_{K'},\C^\ast) \to H^2(K', \Z/p),
$$
which shows that the obstruction to $\varphi=\varphi(K') \, \in \, {\rm Hom}(\Gamma_{K'},\C^\ast)$ being a $p$-th power of another character of $\Gamma_{K'}$ is the class $\partial(\varphi) \in H^2(K', \Z/p)$, where $\partial$ is the connecting morphism from $H^1(\Gamma_{K'},\C^\ast)$ into $H^2(K', \Z/p)$. Put
$$
\beta=\beta(K') : = \, \partial(\varphi(K')), \, \, \, \forall \, \, K' \in \iI.
$$

Note that since the restriction map ${\rm res}_{K/K'}: \, H^i(K',-)\to H^i(K, -)$ commutes with $\partial$, and since $\varphi(K')$ restricts to the trivial character on $\Gamma_K$, the restriction of $\beta(K')$ is trivial in $H^2(K, \Z/p)$. Hence the collection $\{\beta(K') \, \vert \, K'\in \iI\}$ satisfies the hypothesis of Lemma 4.1. Consequently, given any auxiliary finite set $T_0$ of finite places of $F$ which are prime to $p$ and unramified in $K$, with $T$ denoting the set of places of $K$ above $T_0$, we can find a cyclic extensions extension $E$ of $F$ of degree a power of $p$, such that the conclusions (a), (b) and (c) of Lemma 4.1 hold. In particular, for every $K'\in \iI$, and for every subfield $N'$ of
$$
N:= \, KE, \, \, \, {\rm with} \, \, \,  N' \supset K', \, [N:N']=p,
$$
we have
$$
{\rm res}_{N'/K'}(\beta(K')) \, = \, 0.
$$
It is important to note that by construction (cf. Lemma 4.1), the same $E$ works for all $K'$ in $\iI$. Since by definition $\beta(K')=\partial(\varphi)$, and as the restriction map commutes with $\partial$, we get
$$
\partial(\varphi\vert_{\Gamma_{N'}}) \, = \, 0.
$$
Hence $\varphi(K')_{\vert N'}=\psi^p$ for some $\psi\in {\rm Hom}(\Gamma_{N'}, \C^\ast)$, necessarily of order $p^2$ since $\varphi(K')_{\vert N'}$is still of order $p$, cutting out $N$ over $N'$. Put
$$
L:= \, N'(\psi) \, \supset \, N=K'(\varphi)E=N'(\varphi),
$$
where $N'(\nu)$ denotes, for any character $\nu$ of $\Gamma_{N'}$, the cyclic extension of $N'$ cut out by $\nu$. Then $L$ is cyclic of degree $p^2$ over $N'$ and of degree $p$ over $N$. This way we get a collection $\jJ$ of cyclic $p$-extensions $L$ of $N$, one for each $N'$ as above, as $K'$ varies over $\iI$.

For every $\tau$ in $\Gamma_F$, $K^\tau=K$ since $K/F$ is Galois, while $L^\tau$ need not be $L$, though still cyclic of degree $p$ over $N=KE$. Put
$$
\tilde L: = \, \prod_{\tau\in \Gamma_F} \, L^\tau,
$$
and
$$
M: = \, \prod_{L \in \jJ} \, \tilde L.
$$
Then $M/F$ is a finite Galois extension, and since $M$ is a compositum of cyclic $p$-extensions of $K$, $M/K$ is solvable.
So the base changes $\pi_M, \pi'_M$ are defined in $\aA(n,M)$.

\medskip

It suffices to show that for all but a finite number of places $v$ of $K$ which are of degree $p^j$ over $F$, if $u$ is a place of $M$ above $v$, then deg$_{M/F}(u) \geq p^{j+1}$. (The finite number of places $v$ which are ignoring are the ones above $T$.) Pick any such $v$, with $v_1$ denoting the place of $F$ below it. Since $v$ has degree divisible by $p$, its decomposition group (over $F$) contains a cyclic subgroup $H$ of order $p$ in Gal$(K/F)$. Let $K'\in \iI$ correspond to $H$. Let $v_L$ be the place of $L$ below $u$. It suffices to show that $v_L$ has degree $\geq p^{j+1}$ over $F$. If $v_{N'}$, resp. $v'$, is the place of ${N'}$, resp. $K'$, below $v_L$, then the degree of $v_{N'}$ over $F$ is at least as big as the degree of $v'$ over $F$, which is $p^{j-1}$, since $v$ divides $v'$ and is of degree $p$ over $K'$. Putting these together, we see that it suffices to check that $v_L$ has degree $p^2$ over ${N'}$. By construction, ${N'}/L$ is cyclic of degree $p^2$ and moreover, $v_{N'}$ is inert in $KE$, which is cyclic of degree $p$ over ${N'}$. Then $v_{N'}$ must be inert all the way in $L$. Indeed, if it were false, $v_L$ would have degree $p$ over ${N'}$ and its decomposition group in Gal$(L/{N'})$ would be the unique cyclic subgroup $C$, say. Then $KE$ would necessarily be the fixed field of $C$, in which case $v_{N'}$ would split in $KE$, which contradicts the fact that $v_{N'}$ is inert in $KE$. Thus $v_L$ has degree $p^2$ over ${N'}$, and this phenomenon recurs at every possible ${N'}$ as $K'$ varies over $\iI$.

\medskip

This yields part (a) of Proposition 5.1. Part (b) follows as well since by construction, every place of $T$ splits completely in each $L^\tau$, and hence in $M$.

Done.

\medskip

\section{Proof in the general case}

\medskip

Now let $K/F$, $\pi, \pi'$ be as in Thm. A, so that we have intermediate fields $K_j, j\{0,\dots, m\}$ such that
$$
F=K_0 \, \subset \, \dots \, \subset K_m=K,
$$
with each $K_j/K_{j-1}$ Galois ($\forall \, j\geq 1$).

We will now prove Prop. 3.1 in general, which needs to be proved only when
$$
w(\pi, \pi'; K/F) \, > \, 0,
$$
which we will take to be the case.

If $m=1$, $K/F$ is Galois, and the assertion was established in section 5. So let $m>1$ and assume by induction that the assertion holds for $m-1$. Look at the extension $K/K_1$. From the definition of width, we have
$$
w(\pi, \pi'; K/K_1) \, \geq \, w(\pi, \pi'; K/F) \, > \, 0.
$$
Then by induction, for every finite subset $\tilde T$ of $\Sigma_{K/K_1}^1$, prime to the conductor of $\pi, \pi'$ and the discriminant of $K$, there exists a finite extension $N/K$ filtered as $K=N_0 \subset N_1 \subset \dots \subset N_r=N$ with each $N_j/N_{j-1}$ cyclic, such that
\begin{enumerate}
\item[(a)]$w(\pi_N, \pi'_N; N/K_1) \, < \, w(\pi, \pi'; K/K_1))$;
\item[(b)]The places in $\tilde T$ split completely in $N$.
\end{enumerate}

\medskip

We need to conclude the inequality (a) with $K_1$ replaced by $F$. We cannot do this yet, but the offending places $\tilde v$ of $K$ of degree $\leq [(n^2+1)/2]$ are exactly those which are of degree $1$ over $K_1$, but have higher degree over $F$; in other words
$$
\tilde v \mid v, \, \, \, {\rm with} \, \, v\in \cup_{2\leq j\leq [(n^2+1)]/2} \, \Sigma_{K_1/F}^j.\leqno(ast)
$$

Now apply Prop. 5.1 with $K_1=K$ and $T=T_1$, the set of places of $K_1$ below $\tilde T$, to deduce the existence of a solvable extension $M_1/K_1$, which is a succession of cyclic extensions, such that
\begin{enumerate}
\item[(a)]For all but a finite number of $v \in \Sigma_{K/F}^j$ with , we have, for every place $u$ of $M_1$ dividing $v$,
$$
\deg_F(u) \, > \, \deg_F(v);
$$
\item[(b)]The places in $T$ split completely in $M_1$.
\end{enumerate}

\medskip

Consider the compositum $M$ of $K$ and $M_1$, which is a finite solvable extension of $K$ of the type we require. Moreover, by construction, for all but finitely many places $\tilde v$ of degree $1$ over $K_1$ satisfying $(\ast)$, and for all places $\tilde u$ of $M$ above $\tilde v$, one has
$$
\deg_F(\tilde u) \, > \, \deg_F(\tilde v);
$$

We next put
$$
R \, = \, NM,
$$
the compositum of $N$ and $M$.
Evidently,
$$
w(\pi_R, \pi'_R; R/F) \, < \, w(\pi, \pi'; K/F).
$$

Done.

\medskip

\section{Proof of Theorem C}

\medskip

{\it Part (a)}: \, We proceed as in section 2of \cite{Ra2}. We preserve the notations of that paper. It is immediate that it suffices to prove Lemma 2.5 there. s in the proof of that Lemma there, it suffices to show, with $X$ denoting the set of almost all primes $v$ of degree $\geq 2$ where $\pi$ and $\pi'$ are unramified, that
$$
\log L_X(s, \overline \eta \times \eta) \, = \, o(\log\left(\frac{1}{s-1}\right),
$$
for $\eta \in \{\pi, \pi'\}$.

As seen in (2.4) of {\it loc. cit.}, we have, by our choice of $X$,
$$
\log L_X(s, \overline \eta \times \eta) \, = \, \sum_{m=2}^\infty \, c_m m^{-s},
$$
where $c_m=0$ if $m$ is not a power of $Nv$ for some $v \in X$, and when $m$ is of this form, we have
$$
c_m \, = \, \sum_M \, \frac1r \, \sum_{1 \leq i, j \leq n} \, {\overline \alpha_{i,v}}^r{\alpha_{j,v}}^r,
$$
where $M$ is the set of pairs $(v,r)\in X \times \N$ such that $m=Nv^r$, and $\{\alpha_{1,v}, \dots, \alpha_{n,v}\}$ is the Langlands class of $\eta$ at $v$.

By our hypothesis that $\pi, \pi'$ are both semi-tempered, we see that for a positive constant $C_0$,
$$
\vert c_m\vert \, \leq \, C_0\sum_{M=(v,r), {\rm deg}(v)\geq 2} \, Nv^{r(1/2-2\delta)}.
$$
For each $v \in X$, if $\ell$ is the rational prime below it, then as deg$(v)\geq 2$, we have
$$
Nv \, \geq \, \ell^2,
$$
and for each $]ell$ the number of $v$ dividing it is at most the degree of $K$ over $\Q$. It follows that for real $s >1$,
$\log L_X(s, \overline \eta \times \eta)$ is majorized (with $C=C_0[K:\Q]$) by
$$
C\sum_{\ell} \, \ell^{-(2s-1+4\delta)},
$$
which converges at $s=1$. Done.

\medskip
{\it Part (b)}: \, By hypothesis, $\pi$ admits an adjoint transfer $\Pi$ from GL$(n)$ to GL$(n^2-1)$ of $\pi$; it is an isobaric sum of unitary cuspidal automorphic representations of the form
$$
\Pi \, \simeq \, \boxplus_{j=1}^a \, \pi_j,
$$
with each $\Pi_j$ cuspidal on GL$(k_j)/K$ such that $\sum_{j=1}^a \, k_j \, = \, n^2-1$. Applying the Luo-Rudnick-Sarnak bound \cite{LRS} for each $\Pi_j$, we see easily the existence of a constant $t >0$ such that (for almost all $v$)
$$
\vert \beta_{i,v}\vert \, \leq \, Nv^{1/2-t},
$$
where $\{\beta_{1,v}, \dots, \beta_{n^2-1,v}\}$ is the Langlands class of $\Pi_v$. Moreover, since
$$
L(s,\overline \pi \times \pi) \, = \, \zeta_K(s)L(s, \pi, {\rm Ad}),
$$
and since the Langlands class of $\overline \pi_v \boxtimes \pi_v$ is
$\{\overline\alpha_{i,v}\alpha_{j,v} \, \vert \, 1 \leq i,j,\leq n\}$, we deduce that the coordinates of the Langlands class $\{\alpha_{1,v}, \dots, \alpha_{n,v}\}$ of $\pi_v$ satisfy
$$
\vert \alpha_{j,v}\vert \, \leq \, Nv^{1/4 - t/2}, \, \, \forall j \leq n.
$$
Taking $\delta = t/2$ yields the semi-temperedness of $\pi$.

This finishes the proof of Theorem C.

\vskip 0.2in

Dinakar Ramakrishnan

253-37 Caltech

Pasadena, CA 91125, USA.

dinakar@caltech.edu

\end{document}